\documentclass{owrart}

\usepackage{amsmath}
\usepackage{amsfonts}

\newcommand{\RR}{\mathbb R}
\newcommand{\ZZ}{\mathbb Z}
\begin{document}

\begin{talk}[Samuel Fiorini, Gianpaolo Oriolo, and Laura Sanit{\`a}]{Dirk Oliver Theis}
  {The Virtual Private Network Design Problem with Concave Costs}
  {Theis, Dirk Oliver}
  
  \noindent%
  The symmetric Virtual Private Network Design (VPND) problem is concerned with buying capacity on links (edges) in a communication network such that certain
  traffic demands can be met.  The precise definition is below.
  It was shown by Fingerhut, Suri and Turner \cite{FST97} and later, independently, by Gupta, Kleinberg, Kumar, Rastogi and Yener~\cite{GKKRY01} that VPND can
  be solved in polynomial time if it has the so-called tree routing property, that is, each instance has an optimal solution whose support is a tree.  It was
  conjectured that VPND has the tree routing property, see, e.g., Erlebach and R\"uegg \cite{ERL}, Italiano, Leonardi and Oriolo~\cite{ILO02} and Hurkens,
  Keijsper and Stougie~\cite{HKS05}.  The conjecture was recently solved affirmatively by Goyal, Olver and Shepherd~\cite{GOS08} by settling an equivalent
  conjecture, due to Grandoni, Kaibel, Oriolo and Skutella \cite{GKOS08}, claiming that another problem called the Pyramidal Routing (PR) problem has the tree
  routing property.  This fact had previously been established only for cycles \cite{GKOS08,HKS05} and outerplanar graphs \cite{FOST07TR}.

  In recent work we have investigated a natural generalization of VPND where the cost per unit of capacity may decrease if a larger amount of capacity is
  reserved (economies of scale principle).  The growth of the cost of capacity is modelled by a non-decreasing concave function $f$.  We call the problem the
  concave symmetric Virtual Private Network Design (cVPND) problem.  This problem is APX-hard for general $f$ due to the fact that it contains the minimum
  Steiner tree problem as a special case.

  Our main contributions are as follows.  First, as a corollary to Goyal \textit{et al.}~\cite{GOS08}, we show that a generalization of the PR problem which
  we call Concave Routing (CR) problem, and hence also cVPND, has the tree routing property.  Second, we study approximation algorithms for cVPND.  For general
  $f$, using known results on the so-called Single Source Buy at Bulk problem by Grandoni and Italiano~\cite{GI}, we give a randomized $24.92$-approximation
  algorithm.

  \subsection*{Detailed description of the problems}\label{ssec:intro:detailprobs} 

  We now describe the \textit{symmetric Virtual Private Network design} (VPND) problem and its generalization with concave costs, the \textit{concave symmetric
    Virtual Private Network Design} (cVPND) problem.

  The problems have as input a simple, undirected, connected graph $G=(V,E)$ which represents a communication network; a vector $c \in \RR^E_+$ describing the
  edge costs; and a vector $b \in \ZZ^V_+$ of maximum cumulative demands.  A vertex $v$ with $b_v > 0$ is referred to as a \textit{terminal}.  We denote the set
  of terminals by $W$.

  The vertices of $G$ want to communicate with each other.  The exact amount of traffic between pairs of vertices is not known.  However, for each vertex $v$
  the cumulative amount of traffic that $v$ can send or receive is at most its maximum cumulative demand $b_v$.  The aim is to install, at minimum cost, a
  so-called \emph{virtual private network,} as a ``sub-network'' of $G$.  A virtual private network consists of a set of paths $\mathcal P$ containing exactly
  one $u$--$v$ path $P_{uv}$ in $G$ for each unordered pair $\{u,v\}$ of terminals, and a vector $\gamma \in \RR^E_+$ describing the capacity to be installed on
  each edge.  The virtual private network should allow to route any feasible matrix of traffic demands, i.e., whenever $d_{uv}$ is given for each unordered pair
  $\{u,v\}$ of terminals, such that the upper bounds $\sum_{u \in W} d_{uv} \leq b_v$ for the cumulative demands are satisfied, then for each edge $e$, the
  amount of traffic which is routed over edge $e$ is at most the installed capacity:
  \[
  \gamma_e \geq \sum_{\{u,v\} \subseteq W : e\in P_{uv}} d_{uv} \qquad\text{for all edges } e \in E.
  \]

  In the VPND problem, the cost of installing capacity $\gamma_e$ on edge $e$ is $c_e\gamma_e$, while in the cVPND problem, it is $c_ef(\gamma_e)$.  In both
  cases, the cost of a virtual private network is the sum of the installation of capacity on the edges.

  \smallskip%
  We also use what we call the \textit{Concave Routing} (CR) problem, which is a straightforward generalization of the so-called \textit{Pyramidal Routing} (PR)
  problem defined in \cite{GKOS08}.  It, too, is concerned with finding a set of paths satisfying an optimality condition.  We do not define it here.

  \subsection*{Tree routings}

  A feasible solution to one of the problems described above is a \textit{tree solution} if the support of the capacity vector $\gamma$ or, equivalently, the
  union of the paths in $\mathcal P$ induces a tree in $G$.  We say that a problem has the \textit{tree routing property} if, for any instance, there is always
  a tree solution among the optimal solutions.

  In the case of the VPND, the tree routing property, the fact that one can restrict the search for a minimum cost virtual private network to those networks
  with tree support (tree solutions), established in \cite{GOS08} relying on work in \cite{GKOS08}, is the key to the polynomial time algorithm.  In the
  non-linear case, the tree routing property is important in our design of approximation algorithms.  While acknowledging that the tree routing property for
  cVPND can be proven directly from the tree routing property for VPND, we offer an elegant geometric proof for the tree routing property of the Concave Routing
  problem, which then implies the tree routing property of cVPND in much the same way as the tree routing property for PR implies that for VPND.

  \subsection*{Approximation}

  Our approximation algorithms for cVPND rely on the tree routing property.  We reduce the cVPND to an intermediate problem which might be thought of as an
  undirected uncapacitated minimum concave-cost single-source flow problem.

  In the special case, when the function $f$ has the form $f(x) = \min(\mu x, M)$, this flow problem is known as the Single Source Rent or Buy problem, for
  which a randomized $2.92$-approximation algorithm exists, as was shown by Eisenbrand, Grandoni, Rothvo\ss, and Sch\"afer \cite{EGRS}, improving on results by
  Gupta, Kumar, P\`al, and Roughgarden \cite{GKPR}, which can be derandomized with the factor deteriorating to $3.28$.

  For general functions $f$, the situation is a bit more subtle.  Here we use a theorem of Grandoni \& Italiano \cite{GI} an the analysis of their (randomized)
  approximation algorithm for the so-called Single Source Buy at Bulk problem.  This problem is concerned with minimum-cost installations of cables on edges.
  The cables can be selected from a list of cable types, which have characteristic capacities and installation costs.  After breaking the function $f$ up into a
  polynomial sized list of cable types satisfying certain scaling properties, we are able to use said theorem in \cite{GI}, to show that a polynomial time,
  randomized $24.92$-approximation algorithm exists for the afore-mentioned flow problem, and hence also for the cVPND.

\end{talk}
\end{document}